\documentclass[twocolumn,10pt]{article}
\usepackage[dvips]{graphicx}
\begin{document}
\bibliographystyle{prsty}

\title{Destruction of CD4 T Lymphocytes Alone Cannot Account for their Long-term Decrease in AIDS}

\date{\today}
\maketitle

\begin{center}
{\large \bf Yoram Louzoun{\scriptsize (1)},Irun. R. Cohen{\scriptsize (2)}, and Henri Atlan{\scriptsize (3)}}  
\end{center}

\textit{\scriptsize (1) Interdisciplinary Center for Neural Computation, Hebrew University, Jerusalem Israel.  yoraml@alice.nc.huji.ac.il}

\textit{\scriptsize (2) Department of Immunology, The Weizmann Institute of Science, Rehovot, Israel.}

\textit{\scriptsize (3) Human Biology Research Center, Hadassah Hebrew University Hospital, Jerusalem Israel.}

\thispagestyle{empty}
\noindent
\begin{abstract}
Following previous models describing a quasi steady state (QSS) for the evolution of HIV infection and AIDS\cite{obsQSS1,obsQSS2}, we have developed a larger formalism simulating the long-term evolution of the QSS\cite{large_hiv}.  We show that the long-term evolution of AIDS cannot be explained by the destruction alone of CD4 T cells, either directly or indirectly.  The destruction of CD4 T cells can lead only to a QSS with a lower concentration of CD4 T cells, but CD4 destruction cannot generate the sustained long-term decrease in T cells leading to AIDS.  We here suggest some workable explanations.
\end{abstract}

\section{Introduction}
The last few years have witnessed an ongoing debate over the kinetics of AIDS.  Ho et al.\cite{ho} and Wei et al.\cite{wei}  showed that virions and infected CD4 T cells (TL) have a very short life span, and that the HIV and the immune system are in a QSS during the latent stage.  Many mechanisms were proposed to account for the slow decrease in the CD4 TL concentration during the latent stage; most of them involve the direct or the indirect destruction of the CD4 TL by the virus.  For a review and a simplified network representation see Atlan 1996 \cite{atlan_96}.  We here employ a simple mathematical formalism to show that none of these mechanisms can explain the long term decrease in the CD4 TL concentration, and we propose other mechanisms to account for this decrease. 

\section{QSS T Cell Concentration}
In order to describe qualitatively the change in the overall CD4 TL concentration ($T$) during HIV infection, we will use a QSS formalism\cite{QSS_form}.  The QSS formalism assumes that all the mechanisms affecting $T$ can be averaged as a function of $T$ \footnote{This assumption is true for the virions and for the infected CD4 TL  \cite{ho,wei} and for all the cell types reacting directly to the virions or the infected cells.}.
Before infection, $T$ is determined by a production rate ($a$), a proliferation rate ($b$) and a death rate ($\delta_T$). (We will address the possibility of a homeostatic mechanism and saturation in the next section.)
\begin{equation}
\dot{T} = a+(b-\delta_T)T=a-yT;y=-b+\delta_T
\label{eqn:heal1}
\end{equation}

The steady state is $T=a/y$, and the TL concentration reaches this steady state within a time of $O(1/y)$ (Figure 1a). \footnote{The precise solution to this equation is $T={a \over y} +(T(0) - {a\over y})e^{-yt}$.}
During the latent stage of the infection, $T$ is lower than in the healthy state. If the decrease in $T$ QSS is due to a constant destruction of CD4 TL (directly by the virions or by anti CD4 CTL), the equation we get is:

\begin{equation}
\dot{T} = a_1-yT -\gamma T
\label{eqn:dis1}
\end{equation}

Where  $-\gamma T$ represents the destruction rate. The steady state of this equation is $T={a \over y+\gamma}$, and $T$ reaches the steady state in a time $O({1 \over y+\gamma})$, which is faster than $O({ 1\over y})$ \footnote{The precise solution in this case is $T={a \over y+\gamma} +(T(0) - {a\over y+\gamma})e^{-(y+\gamma)t}$.} (Figure 1b). Thus, adding a constant destruction rate can lower the QSS concentration, but this alone cannot produce the long-term decrease in $T$ characteristic of AIDS; $T$ only stabilizes around a new steady state.

One might reason that AIDS could emerge if one assumes that the CD4 TL destruction rate changes as a function of $T$.  In other words imagine that the concentration of the agent (cell or virion) destroying the CD4 TL depends on $T$.  For example, if we have an agent $D$ with a proliferation rate proportional to $T$ and a constant destruction rate, the dynamics of $D$ would be described by: 

\begin{equation}
\dot{D}=xT-\delta_DD
\end{equation}

The QSS concentration of D then would be $D={x \over \delta_D}T$. The evolution of $T$ would then be described as :
\begin{equation}
\dot{T} = a-yT - {x \over \delta_D}T^2
\end{equation}

This case can be generalized to: 
\begin{equation}
\dot{T} = a-yT -\gamma T^n;n>1
\end{equation}
Note that $T$ will still reach a steady state, but more quickly. The steady state value of $T$ will be in this case: $T={a+(n-1)\gamma(a/y)^n \over y+n\gamma(a/y)^{n-1}}$ which converges to the linear steady state for $n=1$. The time at which $T$ reaches a steady state is to first order  $O({1 \over y+\gamma n ({a \over y})^{n-1}})$ (Figure 1c). Thus, neither a constant rate of CD4 TL destruction nor a $T$-dependent rate can explain the long-term fall in TL concentration during HIV infection. This is true even if TL destruction is caused by some indirect mechanism, for example, apoptosis, destruction of CD4 TL by CTL, anergy or other means.

\begin{figure}[ht]
\center
\noindent
\includegraphics[clip, width = 6 cm]{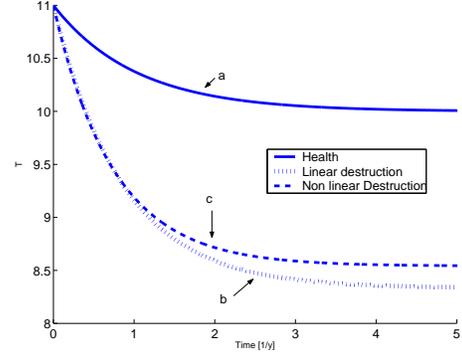}
\caption{The evolution of CD4 TL concentration (a) in health, (b) in the event of a linear destruction rate, and (c) by a non-linear destruction rate.}
\end{figure}

\section{Saturation}
The same reasoning applies even if the CD4 TL steady state is regulated by a homeostatic mechanism.  In this case the proliferation rate can be higher than the death rate, and the total TL population is limited by a higher order mechanism; we thus replace eq \ref{eqn:heal1}  by:

\begin{equation}
\label{eqn:hom1}
\dot{T}=a+yT-\gamma T^2
\end{equation}
This equation can be solved analytically\footnote{The precise solution is: $T(t)={y+\sqrt{y^2+4a\gamma}\tanh{{\gamma \sqrt{y^2+4a\gamma} t -\sqrt{y^2+4a\gamma}C\over 2 \gamma}} \over 2 \gamma}$} and has 2 types of behaviors:

\begin{enumerate}
\item If the T cells are produced (in the thymus or elsewhere) at a constant rate $a > yT$, we can neglect $yT$ and the dynamics of $T$ will be determined by the putative regulatory factor $\gamma T^2$. This factor could be the target of a process that raises the value of $\gamma$ (Figure 2b), or lowers the value of $a$ (Figure 2c); but in either case the rate with which $T$ reaches a steady state is faster than in the healthy state (Figure 2a,2b)\footnote{The steady state is $T=\sqrt{a \over \gamma}$. Starting from an initial value of $T(0)$ the initial evolution is: $T(0) - ( T(0)^2\gamma-a) t$. Thus lowering $a$ or rising  $\gamma$ will both lead to a lower steady state, but this steady state will be aproched faster than in the healthy case. Indeed if ones lower $a$ the evolution rate of $T$ during the last stages of the disease is slower than in the healthy case, but this slow evolution occurs only very close to the steady state, and has no effect in most stages of $T$ evolution}.

\begin{figure}[ht]
\center
\noindent
\includegraphics[clip, width = 6 cm]{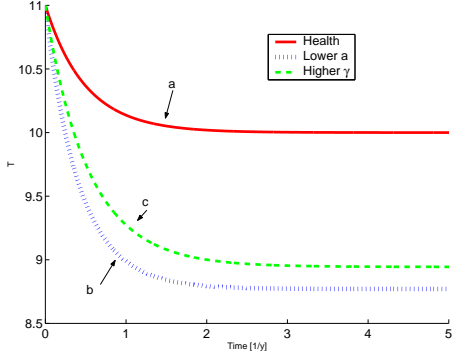}
\caption{The evolution of CD4 TL concentration if it is dominated by a homeostatic mechanism (a) in health, (b) with a stronger homeostatic effect, and (c) with a lower constant source.}
\end{figure}

\item  If, in contrast, the main source of T cells is by way of the proliferation of existing T cells ($a$ is negligible), then the dynamics of $T$ will be dominated by the net proliferation rate $y$. In this case a lower steady state of $T$ can be achieved either by raising $\gamma$, or by lowering $y$ (that is, by adding a constant destruction rate.). As in the cases described above, a new steady state is achieved, but it is achieved faster than in the healthy case (Figure 3 b-e) \footnote{The steady state is $T={y \over \gamma}$. Starting from an initial value of $T(0)$ the initial evolution is: $T(0) -(T(0)\gamma-y) t$. Thus lowering $y$ or rising  $\gamma$ will both lead to a lower steady state, but this steady state will be achieved faster than in the healthy case.}.
\end{enumerate}

To summarize, we can say that any direct or indirect destruction of CD4 TL will lower the value of the TL concentration at the steady state, but such destruction alone will not delay the time it takes to reach this steady state. Thus a destructive mechanism alone cannot account for the very slow decrease in T actually observed in AIDS. 

Moreover, any mechanism describing $T$ as a function of time generates a concave curve. In other words, the computed rate of decrease of T as a function of time decreases. However, the typical observed curve of evolution of AIDS is actually convex (Figure 3f). 

\begin{figure}[ht]
\center
\noindent
\includegraphics[clip, width = 6 cm]{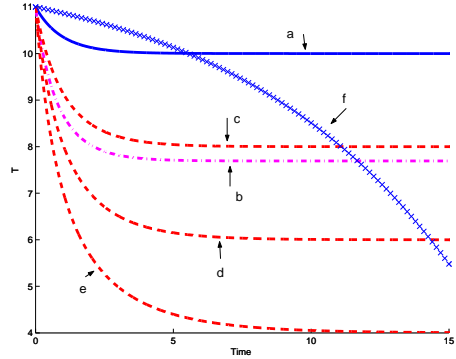}
\caption{The evolution of CD4 TL concentration if it is dominated by proliferation and a homeostatic mechanism (a) in health, (b) with a stronger homeostatic effect, (c-e) with a lower proliferation rate, and (f) the experimentally observed curve} 
\end{figure}

\section{Slow Change in QSS}
The QSS analysis presented in the previous sections shows that any factor we add to the destruction rate of TL will reduce their concentration at the steady state, but such a factor will also reduce the time it takes for $T$ to reach this steady state. Thus, although many mechanisms can explain the decrease in the TL concentration, none of those proposed until now are able to explain why the decrease evolves so slowly. 

The previous analysis fails only if the QSS assumptions are wrong. These assumptions can be wrong only if the process explaining the decrease in TL acts on cells whose lifespan is very long, and only if this process is enhanced as the value of $T$ decreases. 

Thus the CD4 TL concentration can decrease slowly only if all the following conditions are fulfilled (Figure 4):

\begin{itemize}
\item The mechanism of destruction does not act directly on TL.
\item The destruction rate is very slow.
\item The rate of destruction increases when the TL concentration decreases. 
\end{itemize}

These requirements are foreign to most known mechanisms. For example, if TL lose is due to a constant destruction of TL, the number of T cells destroyed will merely decrease when their concentration decreases. Accordingly, the destruction rate can only decrease as the T cell concentration decreases. Thus, for the destruction rate to increase when the concentration of TL decreases, the destruction of TL must feature some positive feedback on itself, and the positive feedback must be stronger than the negative feedback due to the decrease in TL concentration. 

\begin{figure}[ht]
\center
\noindent
\includegraphics[clip, width = 6 cm]{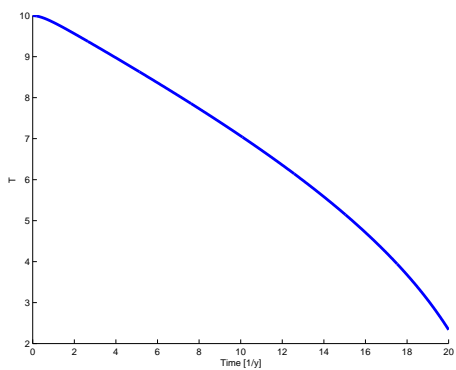}
\caption{The evolution of T as a function of time due to a weak autoimmune destruction of the lymph nodes} 
\end{figure}

\section{Proposed Mechanisms }
We propose here some simple mechanisms that can fulfill the requirements that we showed above to be necessary to explain the long-term evolution of $T$ during HIV infection.  These include:

\begin{enumerate}
\item A slow evolution of the virus toward more destructive variants \cite{urihiv}.
\item A reduction in the CTL proliferation rate by cytokine regulation.
\item A competition between humoral and cellular immunity. 
\item A reduction in the rate by which B cells produce virus-neutralizing antibodies due to destruction of the functional lymph node architecture.
\end{enumerate}

Here are some details:

\begin{itemize}
\item Cytokines - The effect of cytokines on HIV evolution is still controversial \cite{cytshift}.  The observed decrease in TH1 type cytokines and the increase in TH2 type cytokines might be either a marker of the decrease in the number of CTLs or its cause.  The tendency of cytokines to positively influence the proliferation of the T cells that produce them can serve as a positive feedback loop which might balance the effect of the negative feedback caused by the loss of infected cells.Infected T cells secrete TH2 type cytokines and CTLs secrete TH1 type cytokines \cite{TH1cyta}. During the latent stage, the CTL concentration decreases faster than does that of the infected cell population. Hence, the TH1 concentration decreases and the TH2 concentration increases. At some stage, there is an inversion in the cytokine profile and a collapse of CTL production.  This leads to a rapid rise in the viral load and a fast decrease in CD4 due to accelerated infection. 

\item Competition between humoral and cellular immune responses. - The decrease in the concentration of TH1 cytokines and the increase in TH2 cytokine concentration can also be explained by an extrapolation of classical TH1/TH2 competition.  The humoral and cellular immune systems act in parallel against HIV.  Antibodies destroy the virus, while, CTLs destroy the infected cells.  An increase in the B-cell concentration will lead to a decrease in the virus concentration leading, in turn, to a decrease in the concentration of infected cells.  But the CTLs become activated by the infected cells.  Consequently the CTL concentration will decrease when the infected cell concentration decreases.  This mechanism leads to a competition in which raising the humoral response will decrease the cellular response and vice versa.  This simple competition can be enlarged and strengthened if we assume that CTLs, or the CD4 TL that activate CTLs, secrete TH1 type cytokines and that B cells secrete TH2 type cytokines.  Both cytokine types enhance their own production by T helper cells and limit the secretion of the opposite type.  This competition can enlarge an existing competition between B and T cells leading to the collapse of the cellular immune response. 

\item B-Cell depletion - The destruction of lymph nodes and of their germinal centers can limit the proliferation both of B cells and of T cells.  If B-cell antibody production is reduced, the number of free virions can increase despite a reduction in the number of infected cells. The action of such a mechanism might be deduced from the negative correlation between virus and CTL concentration.  The CTL concentration is proportional to the infected cell concentration.  Thus, a negative correlation between CTL and virus concentrations can be due to a negative correlation between the concentrations of virus and infected cells.  The virus concentration rises while the infected cell and CTL concentration is either constant or decreases.  This can be achieved only if the viral lifespan increases.  This increase in lifespan could be due to a reduction in the capacity of B cells to make antibodies that neutralize HIV. 

\end{itemize}

\section{Conclusions}
We have shown that no simple direct or indirect mechanism of TL destruction can explain the long-term decrease in CD4 TL concentration leading to AIDS. The long term decrease in CD4 TL concentration can be explained only by some positive feedback of the destructive mechanism on itself.  We propose four mechanisms, which separately or together, can produce the required positive feedback.  Understanding the precise mechanism that is responsible for the long-term decrease of CD4 TL  is of vital importance.  Treating this mechanism immunologically might help halt the progression to AIDS, in synergy with drug therapy. 

\section{Acknowledgment}
Irun R.Cohen is the incumbent of the Mauerberg Chair in Immunology, Director of the Robrt Koch-Minerva Center for Research in Autoimmune Disease, and Director of the Center for the Study of Emerging Diseases.

Henri Atlan is the incumbent of the Frank Taper chair in medical diagnosis at the Hebrew University in Jerusalem

\end{document}